\documentclass{notices}
\usepackage{amsfonts}%
\usepackage{amssymb}%
\usepackage{graphicx}%
\usepackage{url}%
\usepackage{hyperref}
\usepackage{amsxtra}
\usepackage{amsthm}
\usepackage{amsrefs}
\usepackage{tikz}
\usetikzlibrary{calc}
\usetikzlibrary{patterns}
\newcommand*\rows{20}

\usepackage{float}
\usepackage{wrapfig}

\usepackage{pgfplots}
\pgfplotsset{compat=newest,trig format plots=rad,
    }
\usepgflibrary{fpu}
\usepackage{accents}
\usepackage{enumitem}
\usepackage{subcaption}

\newcommand\Pone{\textrm{P}_{\textrm{I}} }
\newcommand\Ptwo{\textrm{P}_{\textrm{II}} }
\newcommand\Pthree{\textrm{P}_{\textrm{III}}}
\newcommand\Pfour{\textrm{P}_{\textrm{IV}}}
\newcommand\Pfive{\textrm{P}_{\textrm{V}}}
\newcommand\Psix{\textrm{P}_{\textrm{VI}}}

\newcommand{\wup}{\overline w}
\newcommand{\wdo}{\underline w}

\begin{document}

\title{Discrete Painlev\'e Equations}
\author{
        Nalini Joshi
   \affil{
         The author is a professor of mathematics at the University of Sydney.
         Her email address is nalini.joshi@sydney.edu.au.
        }
      }
\date{}
\maketitle

\section{Introduction}
This article is about a class of special functions that cannot be expressed in terms of classical elementary functions. Nevertheless, they arise in a wide range of applications and have surprisingly rich mathematical properties.

In this article, we describe some of these properties, which lie at the intersection of many directions in mathematics, including  dynamical systems theory, differential or difference Galois theory and algebraic geometry. These intersections are expanded in further detail in lectures delivered at UCLA and interested readers may find it helpful to view the videos of the lectures at \url{https://www.math.ucla}\url{.edu/dls/nalini-joshi}\footnote{The author would like to thank Mason Porter and Terry Tao at UCLA for their detailed comments on this article.}. 

\begin{figure}[H] 
\centering
\begin{tikzpicture}[scale=0.36]
            \draw[thin,fill=brown!50] (-0.1,1) rectangle (10.55,3);
            \draw (10,1) -- (10,1.4);
            \node[anchor=north,xshift=0.05 in,yshift=0.04 in,scale=0.5] at (9.9,1.4) {10};
            \draw (0,1) -- (0,1.4);
            \node[anchor=north,xshift=0.04 in,yshift=0.04 in,scale=0.5] at (0,1.4) {0};
            \foreach \x in {1,...,9}{
            \draw (\x,1) -- (\x,1.4) node[anchor=north,xshift=0.04 in,yshift=0.04 in,scale=0.5]{\x};
            }
            \foreach \x in {0.1,0.2,...,10.5}{
            \draw (\x,1) -- (\x,1.075);
            }
            \foreach \x in {0.5,1,...,9.5}{
            \draw (\x,1) -- (\x,1.15);
            }
            \draw (0 in,3) -- (0 in,2.5) node[anchor=north,xshift=0.03 in,yshift=0.03 in,scale=0.52]{0};
            \foreach \x in {1,...,4}{
            \draw (\x in,3) -- (\x in,2.5) node[anchor=north,xshift=0.03 in,yshift=0.03 in,scale=0.52]{\x};
            }
            \foreach \x in {0.1,0.2,...,4.1}{
            \draw (\x in,3) -- (\x in,2.825);
            }
            \foreach \x in {0.5,1,...,3.5}{
            \draw (\x in,3) -- (\x in,2.4);
          };
          \draw[thin,fill=brown!50] (10.55,1) rectangle (21.15,3);
            \draw (20.6,1) -- (20.6,1.4);
            \node[anchor=north,xshift=0.05 in,yshift=0.04 in,scale=0.5] at (20.5,1.4) {10};
            \foreach \x in {0,1,...,9}{
            \draw (10.6+\x,1) -- (10.6+\x,1.4)node[anchor=north,xshift=0.04 in,yshift=0.04 in,scale=0.5]{\x};
            }
            \foreach \x in {0.1,0.2,...,10.5}{
            \draw (10.6+\x,1) -- (10.6+\x,1.075);
            }
            \foreach \x in {0.5,1,...,9.5}{
              \draw (10.6+\x,1) -- (10.6+\x,1.14);
              }
              \draw (4.173228 in,3) -- (4.173228 in,2.5);
              \node[anchor=north,xshift=0.03 in,yshift=0.03 in,scale=0.52] at (4.2 in, 2.5) {0};
            \foreach \x in {1,...,4}{
            \draw (\x in + 4.173228 in,3) -- ( \x in + 4.173228 in,2.5) node[anchor=north,xshift=0.03 in,yshift=0.03 in,scale=0.52]{\x};
            }
          \foreach \x in {0.1,0.2,...,4.1}{
            \draw (\x in + 4.173228 in,3) -- (\x in + 4.173228 in,2.825);
            }
            \foreach \x in {0.5,1,...,3.5}{
              \draw (\x in + 4.173228 in,3) -- (\x in + 4.173228 in,2.4);
              }
         \end{tikzpicture}
\caption{The combined length of two adjacent horizontal rulers is the sum of each individual ruler's length.}\label{f:ruler}
\end{figure}
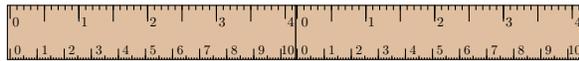

Imagine two rigid rulers, one placed next to the other. (See Figure \ref{f:ruler}.) The length of the combined object is the sum of the lengths of each individual ruler. This is an example of a linear system.

But whether we are observing the heights of colliding waves near a beach,  the ebb and flow of infective cells in a patient, or massive fluctuations in gravitational fields arising from colliding black holes, the corresponding mathematical models turn out not to be linear --- for example, the sum of the heights of travelling waves near a beach before collision is not their combined height at the time of collision. This article is about transcendental functions that solve nonlinear mathematical models arising in such applications.  

The history behind these functions starts, like mathematics, with counting. Addition and subtraction of counting numbers lead to the integers. Multiplication and division lead to rational numbers. Solving polynomial equations with integer coefficients leads to algebraic numbers. But there are numbers called \emph{transcendental numbers} that escape algebraic construction. 
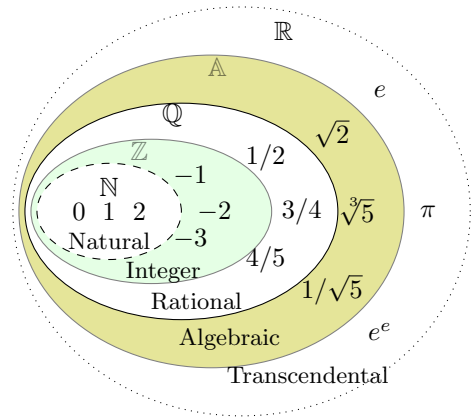
\begin{figure}[H] 
\centering
\begin{tikzpicture}[scale=0.8]
  \draw[dotted] (1.7,0) ellipse (3.8 and 3.4) node at (2.4,3)  {$\mathbb R$};
  \draw[fill=blue!20!yellow,opacity=.5] (1.2,0) ellipse (3.2 and 2.6) node at (1.3,2.4)  {$\mathbb A$};
  \draw[fill=white] (0.7,0) ellipse (2.6 and 1.8) node at (0.55,1.6)  {$\mathbb Q$};
  \draw[fill=green!20!white,opacity=.5] (0.2,0) ellipse (2 and 1.2) node at (0,1.0)  {$\mathbb Z$};
  \draw[fill=white,dashed] (-0.5,0) ellipse (1.2 and 0.8) node at (-0.5,0.4) {$\mathbb N$};
  \node[below] at (-0.5,-0.2) {\small{Natural}} ; 
  \node[below] at (0.4,-0.68) {\small{Integer}} ; 
  \node[below] at (0.95,-1.2) {\small{Rational}} ; 
  \node[below] at (1.5,-1.8) {\small{Algebraic}} ; 
  \node[below] at (2.8,-2.43) {\small{Transcendental}} ;
  \node at (-1,0) {$0$};
  \node at (-0.5,0) {$1$};
  \node at (0,0) {$2$};
  \node at (0.85,0.6) {$-1$};
  \node at (1.25,0) {$-2$};
  \node at (0.85,-0.45) {$-3$};
  \node at (2.1,0.9) {$1/2$};
  \node at (2.7,0) {$3/4$};
  \node at (2.1,-0.8) {$4/5$};
  \node at (3.2,1.3) {$\sqrt 2$};
  \node at (3.6,0) {$\sqrt[\leftroot{-3}\uproot{3}3]{5}$};
  \node at (3.2,-1.3) {$1/\sqrt 5$};
  \node at (4,2) {$e$};
  \node at (4.8,0) {$\pi$};
  \node at (4,-2.0) {$e^e$};
\end{tikzpicture}
\caption{Number systems}\label{f:real}
\end{figure}
Transcendental numbers lie in places on the number line that are like regions in old maps filled with drawings of fabulous monsters. Although the name imbues them with a mystical aura, they arise everywhere in real life, as we know from the presence of the number $\pi$ in almost every part of mathematics. This article is about functions that play a similar role.
\begin{figure}[H] 
\centering
\begin{tikzpicture}[scale=0.84]
  \draw[dotted] (1.45,0) ellipse (3.4 and 3);
  \draw[fill=blue!20!yellow,opacity=.5] (0.7,0) ellipse (2.6 and 1.8);
  \draw[fill=white]  (0.2,0) ellipse (2 and 1.2) ;
  \draw[fill=green!20!white,opacity=.5]  (-0.5,0) ellipse (1.2 and 0.8);
  \node[below] at (-0.5,-0.1) {\small{Polynomial}} ; 
  \node[below] at (0.4,-0.68) {\small{Rational}} ; 
  \node[below] at (0.95,-1.2) {\small{Algebraic}} ; 
  \node[below] at (1.8,-1.8) {\small{Transcendental}} ; 
   \node at (-1.2,0) {$1$};
  \node at (-0.7,0) {$x$};
  \node at (0,0) {$x^2$};
  \node at (0.9,0.8) {$\frac{1}{x^2}$};
  \node at (1.2,0) {$\frac{3+2x}{x}$};
\node at (2.2,1.1) {$\sqrt x$};
  \node at (2.6,0) {$\sqrt[\leftroot{-3}\uproot{3}3]{x}$};
  \node at (2.2,-1.1) {$1/\sqrt x$};
  \node at (3.6,1.7) {$e^x$};
  \node at (4.4,0) {$\Gamma(x)$};
  \node at (3.6,-1.7) {$J_\nu(x)$};
\end{tikzpicture}
\caption{Function classes}\label{f:functions}
\end{figure}
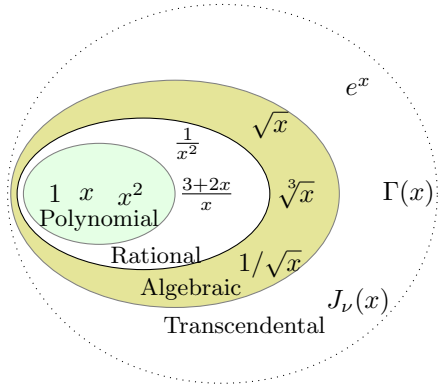

To get to such functions, imagine replacing counting numbers with polynomials.  Multiplication and division of polynomials lead to rational functions. Solving polynomials of more than one variable --- for example, finding solutions $y(x)$ of
\[
a_n(x)y^n+a_{n-1}(x)y^{n-1}+\cdots + a_0(x)=0,
  \]
  where $a_k(x)$ are polynomials of $x$ with integer coefficients --- leads to algebraic functions. Functions that are not algebraic are called \emph{transcendental functions}. 

  Within this class, there are further levels of transcendentality. At the base level sit familiar transcendental functions, such as exponential, trigonometric and logarithmic functions. But there are functions that escape any algebraic operations applied to these or their compositions.

  The exploration of transcendental functions has a long history, reaching back at least to Euler in 1755. One enduring question is whether a newly discovered function can be expressed in terms of earlier known transcendental functions. As for transcendental numbers, the main question here is to ask whether there exist algebraic relations between them. Derivatives or differences of transcendental functions also play a crucial role in answering such questions. 

  Each new level of transcendentality aims to capture functions that escape the previous steps. There is a hierarchy of known cases, which arise as solutions of differential or difference equations. Functions in the hierarchy are categorized according to whether or not they result from operations called \emph{classical operations} applied to functions at a previous step.

  This article focuses on currently known transcendental functions that are furthest away from polynomials, in particular, those at the top of the known hierarchy of functions solving polynomial differential or difference equations. They are solutions of certain nonlinear difference or differential equations, called the {\em discrete Painlev\'e equations} or {\em Painlev\'e equations} \cites{i:56,joshicbms}.

  To describe examples of such equations, we start with basic notation. Difference equations act on the ring $R$ of sequences $(w_n)_{n\in\mathbb N}$ in $\mathbb C$ equipped with an iteration operator $\sigma:R\to R$.  For example,  Euler's Gamma function
\[
\Gamma(z)=\int_0^\infty t^{z-1}e^{-t} dt, 
\]
satisfies a difference equation, given by
$\sigma(w(z))=z\,w(z)$,
where $\sigma(w(z))=w(z+1)$.

Each application of the iteration operator $\sigma$ corresponds to a shift $z\mapsto z+1$ on a line. Given $a\in\mathbb C$, more general shifts  $z\mapsto z+a$ occur on a line and we continue to denote the corresponding iteration by $\sigma(w(z))=w(z+a)$. But, there are more complicated curves on which iterations can be also defined. Two such types of curves and their associated iteration operators have led to extensions of the Gamma function.
One of these is $\sigma_q(w(z))=w(qz)$, with $q\not=0, 1$, $q\in\mathbb C$, which iterates points on a spiral. The other is $\sigma_{\rm ell}$, which iterates points on a curve parametrized by elliptic functions, called an {\em elliptic curve} \cite{cassels1991lmsst}. Given $a, b\in\mathbb C$, elliptic curves have a canonical cubic \begin{figure}[H] 
\begin{subfigure}[b]{0.15\textwidth}
\centering
  \begin{tikzpicture}[scale=0.3]
  \draw[black] (-3,-3) -- (3,3);
 	\filldraw [black] (-2,-2) circle (0.15cm) node[right] {${z-a}$} ;
	\filldraw [black] (0,0) circle (0.15cm) node[right] {${z}$} ;
	\filldraw [black] (2,2) circle (0.15cm) node[right] {${z+a}$} ;
\end{tikzpicture}
\subcaption{$\sigma$}\label{f:iter1}
\end{subfigure}
\begin{subfigure}[b]{0.15\textwidth}
\centering
  \begin{tikzpicture}[scale=0.3]
 \draw[black,smooth,variable=\x]   plot [domain=-3.2:1.2] ( {2*\x*cos(\x r)}, {2*\x*sin(\x r)} )  ;
\filldraw [black] ( {-6*cos(-3 r)}, {-6*sin(-3 r)} )   circle (0.15cm)  node[left] {$q z$};
\filldraw [black] ( {-4*cos(-4 r)}, {4.65*sin(-4 r)} )   circle (0.15cm) node[above] {$z$};
\filldraw [black] ( {-2.4*cos(-2 r)}, {-1.8*sin(-2 r)} )   circle (0.15cm) node[right] {$z/q$} ;
\end{tikzpicture}
\subcaption{$\sigma_q$}\label{f:iter2}
\end{subfigure}
\begin{subfigure}[b]{0.15\textwidth}
\centering
  \begin{tikzpicture}[scale=0.3]
 \draw[black,smooth,variable=\x]   plot [domain=-3.1036:1.3] (\x, {sqrt( ((\x+1)^3-3*(\x+1)+3 ) ) })  ;
 \draw[black,smooth,variable=\x]   plot [domain=-3.1036:1.3] (\x, {-sqrt( ((\x+1)^3-3*(\x+1)+3 ) ) })  ;
\filldraw [black] ( -3, 1 )   circle (0.15cm) node[anchor=south east] {$P_0$} ;
\filldraw [black] ( -0.84,-1.5887 )   circle (0.15cm)  node[anchor=north west] {$P_2$} ;\filldraw [black] ( 1, 2.236 )   circle (0.15cm)  node[right] {$P_1$} ;
\draw[black,dashed] ( -3.8, 0.7528 ) -- ( 1.8, 2.4832 );
\draw[black,dashed] ( -0.84,-2 ) -- ( -0.84,2 );

\end{tikzpicture}
\subcaption{$\sigma_{\rm ell}$}\label{f:iter3}
\end{subfigure}
\caption{Three types of curves, a line, a spiral, and an elliptic curve,  on which iterations are well defined.}\label{f:iter}
\end{figure}
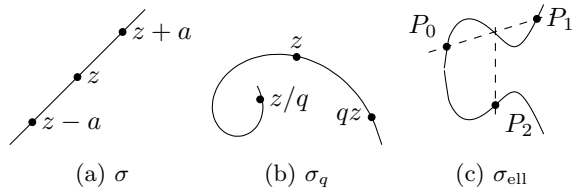\noindent  form $y^2=4x^3-a x - b$. Assume the curve is non-singular, i.e., $a^3\not=27 b^2$. On such a curve, the iteration proceeds as follows: start with two points $P_0$, $P_1$ on the curve. Then the line passing through them intersects the curve again, producing a third point. The iteration of $\sigma_{\rm ell}$ is the the reflection of this third point across the horizontal axis, given by $P_2$ in Figure \ref{f:iter}(\subref{f:iter3}). (This action is also called the addition formula or group action on the curve.) Difference equations with iteration operator $\sigma_q$ are called {\em $q$-difference equations}, while those with iteration operator $\sigma_{\rm ell}$ are called {\em elliptic difference equations}.

All three iteration operators occur in the list of discrete Painlev\'e equations. There are 22 classes of such equations described by Sakai \cite{s:01} through an algebro-geometric study of dynamical systems. A concise list of canonical cases can be found in \cite{joshicbms}*{Appendix D}. Discrete Painlev\'e equations have striking mathematical properties. We illustrate them here for selected examples, which are listed below.

For conciseness, we write $\sigma(w)=\wup$, $\sigma^{-1}(w)=\wdo$, for all three types of iterations. 
  \begingroup
  \allowdisplaybreaks
  \begin{align}
 \label{eq:dp1}&{\rm d}\Pone:\ w\,(\wup+w+\wdo)=a\,n+b+c\,w,\\
\label{eq:qp1}&{\rm q}\Pone :\quad  \overline w\,\underline w =\frac{1}{w}-\,\frac{1}{a\;q^n\, w^2},\\
\nonumber &{\rm e}\textrm{RCG}:\ \\
\nonumber &{\rm cn}(\gamma_n){\rm dn}(\gamma_n)\big(1-k^2{\rm sn}^4(z_n)\big)w_n\big(w_{n+1}+w_{n-1}\big)\\
\nonumber &-{\rm cn}(z_n){\rm dn}(z_n)\big(1-k^2{\rm sn}^2(z_n){\rm sn}^2(\gamma_n)\big)\cdot\\
\nonumber &\phantom{-{\rm cn}(z_n){\rm dn}(z_n)\big(1}\cdot\big(w_{n+1}w_{n-1}+{w_n}^2\big)\\
\nonumber &+\big({\rm cn}^2(z_n)-{\rm cn}^2(\gamma_n)\big){\rm cn}(z_n){\rm dn}(z_n)\cdot\\
&\phantom{-{\rm cn}(z_n){\rm dn}(z_n)\big(1}\cdot\big(1+k^2{w_n}^2w_{n+1}w_{n-1}\big)\,=0,
\end{align}
\endgroup
where
\[
   z_n=(\gamma_e+\gamma_o)n+\omega , \,\,\,\,\,
         \gamma_n=\begin{cases}\gamma_e ,\ \textrm{for}\ n=2j,\\
           \gamma_o ,\ \textrm{for}\ n=2j+1,
           \end{cases}
         \]
  with $\gamma_e$, $\gamma_o$, $a$, $b$, $c$, $d$, $q$ being constants, with $q\not=0,
  1$ and where cn, dn, sn denote Jacobi elliptic functions, i.e., doubly periodic, meromorphic functions that parametrize biquadratic curves.   

  Consider, for example, the 
autonomous case of q$\Pone$, which arises in the limit $|q^n|\to\infty$:
\begin{equation}\label{eq: qp1 auto}
 \overline{w}\,{w}\,\underline{w} =1.
\end{equation}\begin{figure}[H]
\centering
\includegraphics[height=60mm]{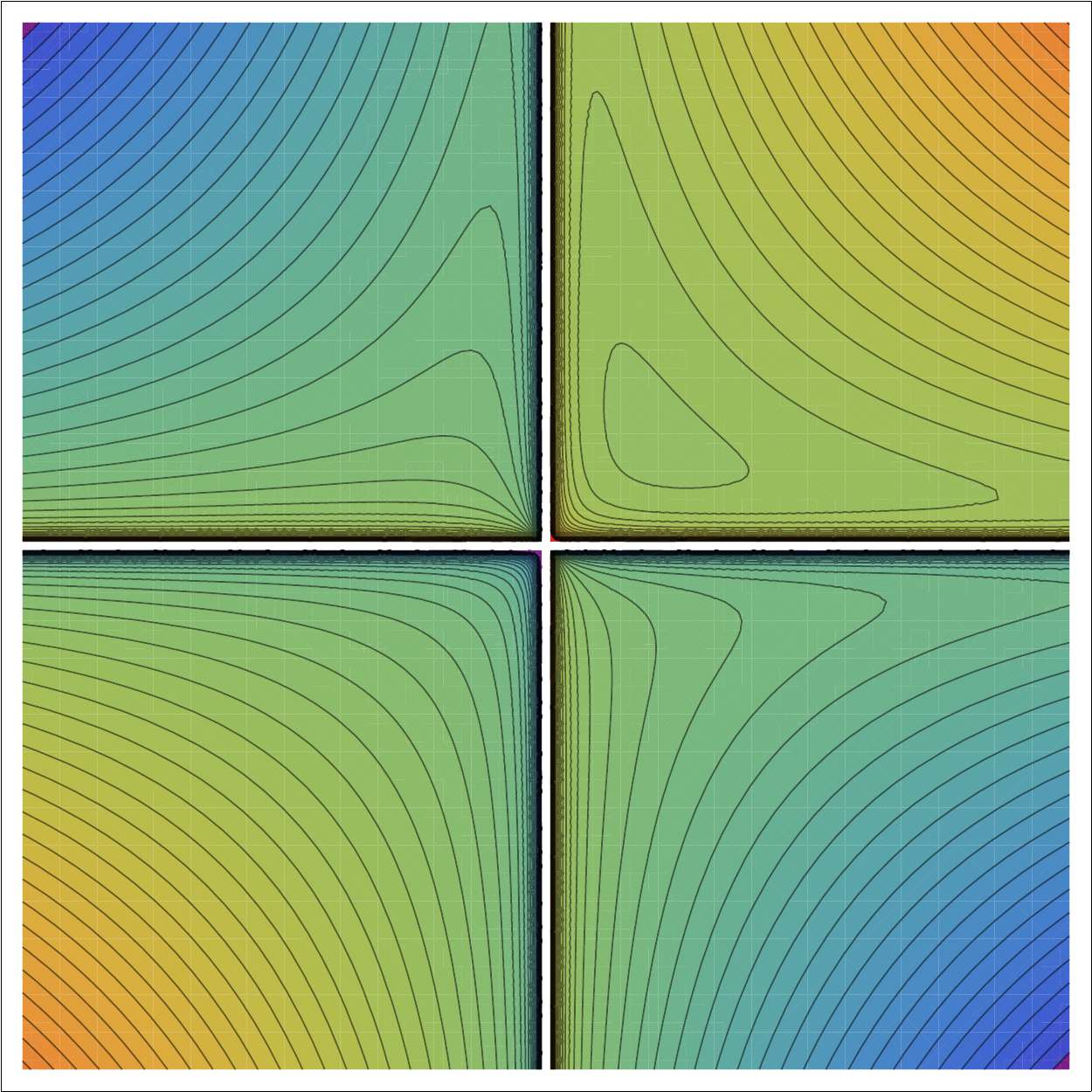}
\caption{Real contour lines of the function $K(x,y)$ defined in Equation \eqref{eq: ham auto}. An initial value $(x,y)=(w_0, w_1)$ lies on one of these curves and the corresponding solution of Equation \eqref{eq: qp1 auto} is iterated on that curve.}\label{fig:ap1curves}
\end{figure}
This equation has an invariant defined by
\begin{equation}\label{eq: ham auto}
K(x, y) =\frac{x^2\,y^2+x+y}{x\,y},
\end{equation}
which satisfies $K(\overline{w}, w)-K({w}, \underline{w})=0$ when $w$ satisfies Equation \eqref{eq: qp1 auto}. Given an initial value $(w_0, w_1)$, the solution orbit, or trajectory, $\{w_n\}_{n=0}^\infty$ lies on the curve defined by $K(x,y)=\kappa:=K(w_0,w_1)$.
This leads to a one-parameter family of curves defined as the zero set of
\begin{equation}\label{eq:ap1curves}
f(x, y)=x^2\,y^2+x+y-\kappa\,xy,
  \end{equation}
which is called a \emph{pencil} of curves  \cite{griffiths1989introduction}. Each curve in this pencil can be mapped to an elliptic curve \cite{cassels1991lmsst}. A two-dimensional real slice of the pencil is illustrated in Figure \ref{fig:ap1curves}.

\begin{figure}[H]
\centering
\includegraphics[height=70mm]{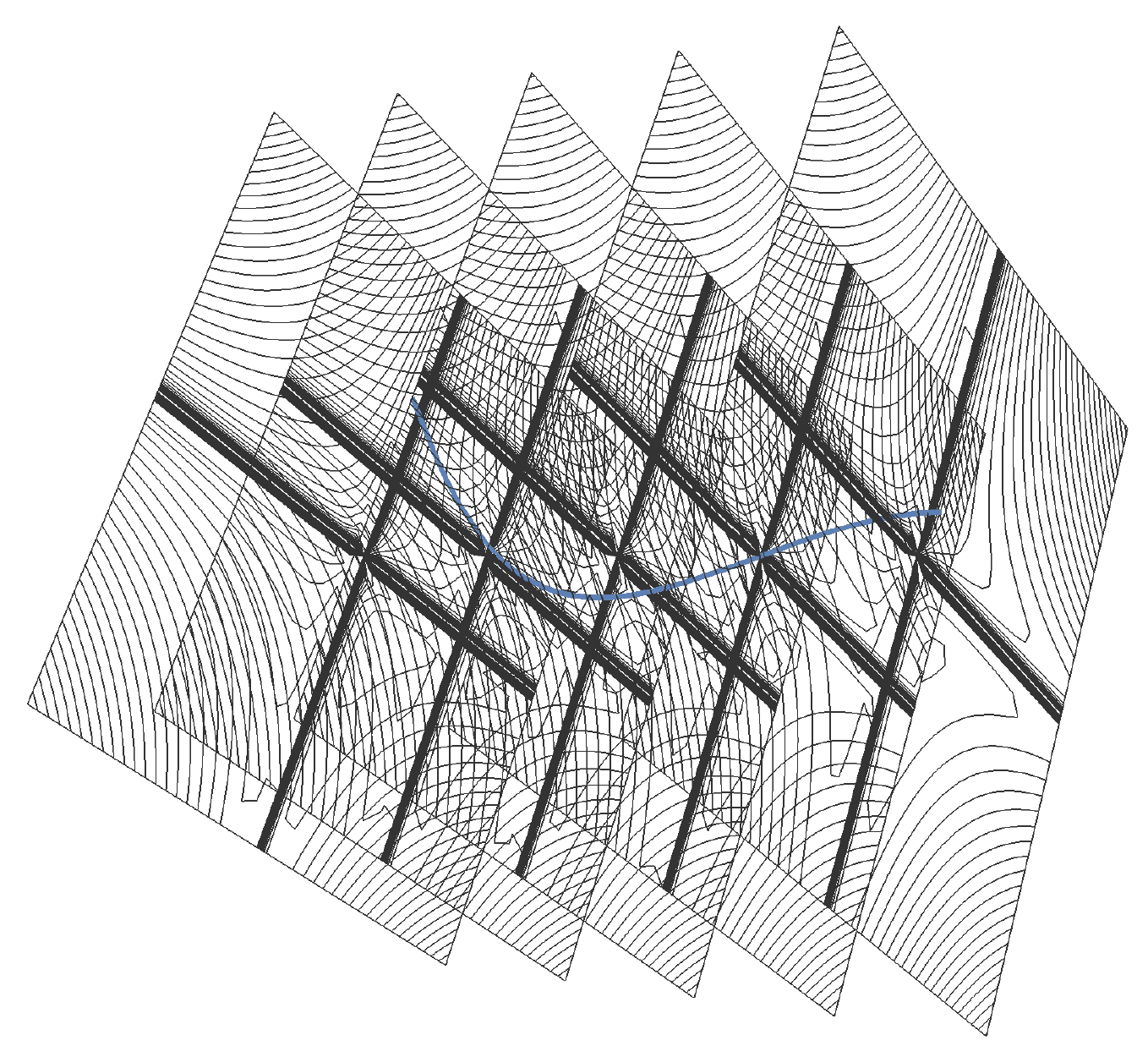}
\caption{A real solution trajectory of q$\Pone$ threading through contour lines of $K$ defined in Equation \eqref{eq: ham auto}.}\label{fig:qp1curves}
\end{figure}

  The function $K$ can be defined by Equation \eqref{eq: ham auto} even when $z=q^n$ is not at infinity, but it is no longer invariant on solution trajectories. Instead, we have
  \[
K(\overline{w}, w)-K({w}, \underline{w})= -\,\frac{1}{z}\frac{w(\overline w - \underline w)}{\bigl(w-1/z\bigr)},
\]
and so for arbitrarily large $|z|$, and bounded values of $w$, $\overline w$, $\underline w$,  $K$ is slowly varying with $z$. A solution trajectory of q$\Pone$ then moves from a point on a contour shown in Figure \ref{fig:ap1curves} to a point on another contour as $n$ is iterated. Such trajectories are like threads that link one contour diagram to another, as shown in Figure \ref{fig:qp1curves}.
Figure \ref{fig:qp1curves} is a local (real) snapshot of a well known structure in mathematics called a {\em foliated vector bundle}. Each slanted contour diagram in the figure is a {\em fibre} of a vector bundle  \cite{milnor} and it fits alongside others locally like parallel planes. Solution trajectories that intersect each fibre transversely form a {\em foliation} of the vector bundle \cite{milnorstasheff}.

In the case of discrete Painlev\'e equations, each fibre of this foliated vector bundle is an elliptic surface, because the invariant contours of the corresponding autonomous system are elliptic curves. Given an initial value, the corresponding solution follows a trajectory that pierces each fibre.  While such trajectories may be locally well defined, this is not guaranteed everywhere. 

We illustrate this for the example of  q$\Pone$, which in system form\footnote{We have only written the forward iteration here for conciseness, but to be complete we need to consider both forward and backward iterations of $(x, y)$.} is given by
\begin{equation}\label{eq:qp1system}
  \begin{pmatrix}\overline x\\ \overline y\end{pmatrix}=
  \begin{pmatrix}\dfrac{z x-1}{z x^2y}\\ x\end{pmatrix},
\end{equation}
where $x=w$, $y=\underline w$, and $z=aq^n$. When $x$ or $y$ are arbitrarily small, their iterates become unbounded, but this can be handled by embedding the system into two-dimensional projective space $\mathbb P^1\times \mathbb P^1$. However, this is not sufficient to handle all problematic initial values, such as the point $b_0:(x, y)=(1/z,0)$ where the right side of Equation \eqref{eq:qp1system} becomes undefined.

For the autonomous case \eqref{eq: qp1 auto}, this is a {\em base point}, i.e., a point where all curves of the elliptic pencil \eqref{eq:ap1curves} intersect. An initial point on each such curve is transported along it with constant speed, like a car travelling on a road. But because all the curves intersect at the origin, there is a gridlock or traffic jam that prevents the cars on each line from moving through this point. The operation of {\em blow-up} or {\em resolution}, which was known to Newton \cite{abhyankar1976historical}, redefines the roads through that impasse, so that each car can travel continuously through it on its journey. See Section \ref{s:ivs} for further detail.

When all base points are resolved, the resulting space is called a space of {\em initial values}. There is a way to characterize this space uniquely for each Painlev\'e and discrete Painlev\'e equation. To do so, we have to understand the self-intersection numbers of lines. Each time a point is blown up on a line, the self-intersection number of that line decreases by 1. For example, a vertical line, denoted by $H_x$ (for constant $x$-value) does not intersect itself in $\mathbb P^1\times \mathbb P^1$. However, when a base point lying on $H_x$ is blown up, the self-intersection number of this line becomes $-1$.  Two separate blow-ups on $H_x$ would lead to a self-intersection number of $-2$.

The resulting resolved initial-value space can be characterised by a reflection group, through an isomorphism known as {\em Duval correspondence} or {\em McKay correspondence}. Here lines of self-intersection number $-2$ play a special role -- they are mapped to simple roots of reflection groups, with their intersections being mapped to inner products.  An example of reflection groups is given in Section \ref{s:arg} and shown to give rise to d$\Pone$ through a translation operator on the lattice generated by such a group.

There are many more properties than we can cover in this article. Each discrete Painlev\'e equation is a compatibility condition for a pair of associated linear problems (called Lax pairs) and Riemann--Hilbert-type methods applied to these linear problems provide information about the solutions.  And, as remarked above, their generic solutions are ``more trancendental'' than any previously known special functions. (See Section \ref{s:ladder}.) In the following, we describe some of these distinctive properties, using the examples listed above for illustration. Further properties and examples are discussed in detail in \cite{joshicbms}.

\section{The beginnings}\label{s:ori}
In 1939, Shohat \cite{s:39} discovered a curious relationship between Equation \eqref{eq:dp1} and orthogonal polynomials. 

Consider the monic Hermite polynomials $\Phi_n(x)$, where $\Phi_0(x)=1$, $\Phi_1(x)=x$, $\Phi_2(x)=x^2-1$, $\ldots$ From the orthogonality relations
\[
\int_{-\infty}^\infty \Phi_n(x)\Phi_m(x) e^{-x^2/2}dx=\sqrt{2\pi}n!\delta_{nm},
\]
where $\delta_{nm}$ is the Kronecker delta, we can deduce a difference equation for $\Phi_n(x)$ (commonly known as a 3-term recurrence relation)
\[
\Phi_{n+1}(x)-x\,\Phi_n(x)+n\Phi_{n-1}=0.
\]

Shohat extended this class to include weight functions of the form
$p(x)=\exp(-x^4/4)$ on $\mathbb R$ and obtained the recurrence relation
\begin{equation}\label{eq:Shohat3term}
  \Phi_n(x)-(x-c_n)\,\Phi_{n-1}(x)+\lambda_n\,\Phi_{n-2}(x)=0, 
\end{equation}
for $n\ge 2$, where $c_n$ are independent of $x$, and deduced a nonlinear difference equation for $\lambda_n$. Extending the weight function further to $p(x)=\exp(-x^4/4+tx^2)$, this equation for $\lambda_n$ becomes
\begin{equation}\label{eq:ShohatdP1}
  \lambda_{n+2}\bigl(\lambda_{n+1}+\lambda_{n+2}+\lambda_{n+3}\bigr)=n+2t\lambda_n.
\end{equation}
This difference equation is a special case of the discrete Painlev\'e equation \eqref{eq:dp1} and is called the \emph{string equation} in physics, due to its appearance in the Hermitian random matrix model of quantum gravity. This is one instance of many applications related to random matrix theory in which solutions of the Painlev\'e and discrete Painlev\'e equations play critical roles \cite{forrester}. 

Denote the solution of Equation \eqref{eq:ShohatdP1} as $\lambda_n(t)$. We can also deduce a differential equation for $\lambda_n(t)$ as a function of $t$. It turns out to be one of the classical Painlev\'e equations \cite{i:56}, which had an entirely different beginning.

In the late \ 19$^{\rm th}$ century,  Picard, Painlev\'e and other mathematicians were engaged in a search for new functions with properties that generalized those of elliptic functions. Their results led to six canonical classes of second-order nonlinear ODEs now called the Painlev\'e equations:
 \begingroup
  \allowdisplaybreaks
  \begin{align*}
\Pone:\quad w''&=6w^2+t\,,\\
\Ptwo:\quad w''&=2w^3+tw+\alpha\,,\\
\Pthree:\quad w''&=\frac{{w'}^2}{w}-\frac{w}{t}+\frac{1}{t}(\alpha w^2+\beta)
      +\gamma w^3+\frac{\delta}{w}\,,\\
\Pfour:\quad w''&=\frac{{w'}^2}{2w}+\frac{3 w^3}{2}+4tw^2+ 2(t^2-\alpha)w-\,\frac{\beta^2}{2w}\,,\\
    \Pfive:\quad w''&=\Biggl(\frac{1}{2w}+\frac{1}{w-1}\Biggr){w'}^2\\
    &\quad -\frac{w'}{t}
      +\frac{(w-1)^2}{t^2 w}(\alpha w^2+\beta)
         +\frac{\gamma w}{t} \\
         &\quad\quad+\frac{\delta w(w+1)}{w-1}\,, \\
\Psix:\quad w''&=\frac{1}{2}\Biggl(\frac{1}{w}+\frac{1}{w-1}+\frac{1}{w-t}\Biggr){w'}^2\\
&\quad -\Biggl(\frac{1}{t}+\frac{1}{t-1}+\frac{1}{w-t}\Biggr)w'\\
               &\quad+\frac{w(w-1)(w-t)}{t^2(t-1)^2}\times\\
    &\quad \times\Biggl(\alpha+ \frac{\beta t}{w^2}
         +\frac{\gamma (t-1)}{(w-1)^2}+\frac{\delta t(t-1)}{ (w-t)^2}\Biggr)\,,
\end{align*}
\endgroup
where $\alpha$, $\beta$, $\gamma$, $\delta$ are constants and $w$ is a function of $t$, with primes denoting derivatives with respect to $t$.

The first difficulty in the search mounted by Painlev\'e was to overcome possible obstructions to analytic continuability created by movable singularities. To illustrate, consider
\[
y'+y^2=0\,,
\]
which has explicit solutions $y(x)=1/(x-x_0)$, for arbitrary $x_0$. (It also has an indentically zero solution.) Such poles are said to be {\em movable}, because their locations are not fixed by the coefficients of the differential equations they satisfy, but instead they move as initial conditions change. To allow analytic continuation of all locally defined solutions,  Painlev\'e imposed a necessary condition that all solutions should be single-valued around all movable singularities. In the modern literature, this condition has been restricted further to require that all movable singularities should be poles and the restricted definition is now widely called the \emph{Painlev\'e property} \cite{joshicbms}.

Painlev\'e devised an ingenious idea to drive the search for ODEs that satisfy his eponymous property.  This idea is based on creating a transformation of variables for a given ODE in the class
\begin{equation}\label{eqn:Pclass}
y''=F(y',y,t),
\end{equation}
where $F$ is rational in $y'$ and $y$ and analytic in a domain in $t$. The transformation incorporates a parameter\footnote{This is not necessarily a parameter in the Painlev\'e equations listed above.} $\alpha\in\mathbb C$ analytically in a domain around $\alpha=0$. Necessary conditions then arise from investigating the limiting case $|\alpha|\to0$ and identifying those cases that fail to have the Painlev\'e property. Because the transformed equation is analytic in $\alpha$, this necessarily implies that the full equation \eqref{eqn:Pclass} also fails to have the required property.

But Painlev\'e's first attempt at classifying ODEs had gaps that were identified and later filled by his student Gambier. Gambier also pointed out that an ODE identified independently by R. Fuchs was the most general case possible in the classification. Painlev\'e also gave sufficient conditions to prove that solutions of the ODEs in each class are globally analytically continuable.\footnote{There are now differing opinions about the validity of his method of proof \cite[\S 1.2]{joshicbms}.}

The proof that generic solutions are ``new'' transcendental functions is non-trivial. An essential step is to show that they cannot be expressed as solutions of algebraic equations with previously known functions as coefficients. Painlev\'e tackled this step by showing that a solution $y$ and its derivative $y'$ cannot satisfy an algebraic equation. We explain these arguments in Section \ref{s:ladder}. But further steps needed an extension of Galois theory for number fields to function fields involving derivatives or differences. There are now several proofs; see \cite{umemura2007invitation}.

The geometric study of initial value spaces of the Painlev\'e equations began in the 1970s, in the French school led by Raymond G\'erard in Strasbourg. This culminated in Okamoto's construction of the compactified and resolved spaces of initial values \cite{okamoto1979} of the six Painlev\'e equations. In 1992, researchers recognized a difference equation arising in a model of quantum gravity as a discrete Painlev\'e equation. This stimulated developments \cite{joshicbms} leading to Sakai's geometric approach \cite{s:01} resulting in 22 classes of discrete Painlev\'e equations.

\section{Affine reflection groups}\label{s:arg}
In this section, we describe an algebraic view of discrete Painlev\'e equations that provides a useful perspective of the structure and dynamics of their solutions. For any given differential (or difference) equation, an enduring question in mathematics is to search for its symmetries -- i.e., transformations under which the equations remain invariant. For the discrete Painlev\'e equations, symmetries are given by affine reflection groups that are closely related to the structure of their initial value spaces. 

When a differential (or difference) equation includes parameters, a symmetry has the following realization. Given a solution corresponding to one choice of parameters, can it be transformed to a solution of another copy of the equation with a different choice of parameters? We illustrate this idea here for $\Pfour$.

\begin{figure}[H]
  \centering
  \begin{tikzpicture}
    \begin{scope}
    \clip (2.6,1.5) rectangle (7.4cm,5.4cm);
          \filldraw[color=gray]  (4.02,3.17) circle (0.05);
          \draw[-latex,thick]  (4.02,3.17) -- (5.01,3.17);
          \node[above] (T) at (4.52,3.17) {$T$};
           \filldraw[color=gray]  (5.02,3.17) circle (0.05);
   \foreach \row in {0, 1, ...,\rows} {
        \draw ($\row*(0.5, {0.5*sqrt(3)})$) -- ($(\rows,0)+\row*(-0.5, {0.5*sqrt(3)})$);
        \draw ($\row*(1, 0)$) -- ($(\rows/2,{\rows/2*sqrt(3)})+\row*(0.5,{-0.5*sqrt(3)})$);
        \draw ($\row*(1, 0)$) -- ($(0,0)+\row*(0.5,{0.5*sqrt(3)})$);
      }
      \end{scope}
    \end{tikzpicture}
    \caption{The arrow in the figure denotes the result of a translation $T$ on a point (drawn in gray) in the triangular ($A_2^{(1)}$) lattice. $T$ shifts each point in a triangle along a distinguished direction (parallel to an edge of the triangle) by a fixed distance (given by the spacing between two parallel lines).}\label{fig:a2}
  \end{figure}
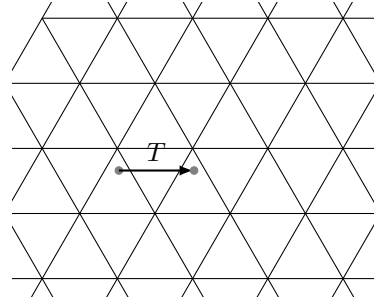

Let $\{\alpha_n\}$ and $\{\beta_n\}$ denote sequences of parameters and $w_n$ the corresponding solution of $\Pfour$. Classical results show that there exist transformations:
\begin{equation*}
\Pfour(w_n, t; \alpha_n, \beta_n)\mapsto \Pfour(w_{n\pm 1}, t; \alpha_{n\pm 1}, \beta_{n\pm 1})\,,
\end{equation*}
where
\begin{subequations}\label{eq:abP4}
  \begin{align}
    \alpha_n&=-\,\frac{n}{2}-\,\frac{c_0}{2}+\frac32 c_1(-1)^n\,,\\
    \beta_n&=n+c_0+ c_1(-1)^n\,,
  \end{align}
\end{subequations}
$c_0$, $c_1$ are given constants, and the solutions are related by
\begin{subequations}\label{eq:bt4}
  \begin{align}
  \label{eq:bt4+} 2 w_n w_{n+1}&=-\,w_n'-\,w_n^2-2tw_n+\beta_n\,,\\
   2 w_n w_{n-1}&=w_n'-\,w_n^2-2tw_n+\beta_n\,.
  \end{align}
\end{subequations}
Adding the two equations \eqref{eq:bt4}, we get
\[
 2 w_n \bigl(w_{n+1}+w_n+w_{n-1}\bigr)=-4t w_n+2 \beta_n,
\]
a generalization of d$\Pone$ with $c=-2t$, $a=1$, and $b=c_0$. 
Okamoto \cite{okamoto1986studiesIII} showed that such transformations form an affine reflection group of type $A_2^{(1)}$, described below. Moreover, the iteration $n\mapsto n\pm 1$ in Equations \eqref{eq:bt4} correspond to translations on the  $A_2^{(1)}$ root lattice, as indicated in Figure \ref{fig:a2}.

To see how affine reflection groups arise, consider a symmetric system equivalent to $\Pfour$ \cite{ny:99}:
\begin{equation}\label{eq:nyp4}
  \begin{aligned}
f_0'&=f_0(f_1-f_2)+\gamma_0\,,\\
f_1'&=f_1(f_2-f_0)+\gamma_1\,,\\
f_2'&=f_2(f_0-f_1)+\gamma_2\,,
\end{aligned}
\end{equation}
where $f_0+f_1+f_2=t$, $\gamma_0+\gamma_1+\gamma_2=1$, and primes 
\noindent
indicate differentiation in $t$. By eliminating $f_0$, $f_2$, it can be shown that $\Pfour$ holds for $f_1$ with $f_1=-w/\sqrt 2$, $t\mapsto \sqrt 2t$, $\alpha=\gamma_0-\gamma_2$ and $\beta=\pm\gamma_1$.

Define operations $s_0$, $s_1$, $s_2$ and $\pi$ on Equations \eqref{eq:nyp4} by
\begin{align*}
  s_i(\gamma_i)&=-\gamma_i\,,&& s_i(\gamma_j)=\gamma_j+\gamma_i\,, \ j=i\pm 1\,, \\
  s_i(f_i)&=f_i\,, && s_i(f_j)=f_j\pm\frac{\gamma_i}{f_i}\,, \ j=i\pm 1\,, \\
 \pi(\gamma_j)&=\gamma_{j+1}\,, && \pi(f_j)=f_{j+1}\,,
\end{align*}
for $i, j\in\mathbb N\ ({\rm mod}3)$. The operators $s_i$ satisfy the {\em Coxeter relations}
\[
s_i^2=1\,,\quad (s_is_{i+1})^3=1\,,\quad i\in \mathbb Z/3\,\mathbb Z ,
\]
which generate the affine reflection group $W=\langle s_0, s_1, s_2\rangle =A_2^{(1)}$. Each $s_j$, $j\in \mathbb Z/3\,\mathbb Z$, represents a reflection across a line in the lattice shown in Figure \ref{fig:a2}. The operator $\pi$, called a diagram automorphism, acts in the following way:
\[\pi^3=1\,,\quad \pi s_i = s_{i+1}\pi\,,\quad i=0, 1, 2\,.\]
Augmenting $W$ by $\pi$ leads to the extended group \[\widetilde W=\langle s_0, s_1, s_2, \pi\rangle \,.\]
These transformations are collectively called {\em B\"acklund transformations} of $\Pfour$.

Given a triangle in the triangular lattice, there are three fundamental lines $\ell_0, \ell_1, \ell_2$ given by extending its edges. Given a point in the lattice, let $\alpha_j$ be given by the orthogonal distance from that point to $\ell_j$, $j=0, 1, 2$. This leads to a natural coordinate system $(\alpha_0, \alpha_1, \alpha_2)$ for points on the lattice.

The lines $\ell_0, \ell_1, \ell_2$ also provide directions along which translations act. So there are three translations on the triangular lattice in Figure \ref{fig:a2}. One of these is given by $T=\pi\,s_2\,s_1$. The remaining translations are equivalent to this one under conjugation by group operations.

The actions of $s_i$, $i=0, 1, 2$, and $\pi$ on the coordinates $(\alpha_0, \alpha_1, \alpha_2)$ can be found explicitly by using Euclidean geometry.  We find
\[
T(\alpha_0)=\alpha_0+1\,,\ T(\alpha_1)=\alpha_1-1\,,\ T(\alpha_2)=\alpha_2\,,
\]
and
\[
  \begin{aligned}
    T(f_1)&=t-f_0-f_1-\,\frac{\alpha_0}{f_0}\,,\\
    T^{-1}(f_0)&=t-f_0-f_1+\,\frac{\alpha_1}{f_1}\,.\\    
  \end{aligned}
\]
Setting $x_n=T^n(f_1)$, $y_n=T^n(f_0)$, we obtain a system of difference equations:
\begin{align}
x_{n+1}&=t-y_n-x_n-\frac{\alpha_0+n}{y_n}\,,\\
y_{n-1}&=t-y_n-x_n+\frac{\alpha_1+n}{x_n}\,,
\end{align}
which contains d$\Pone$. This result illustrates a crucial point in the construction of discrete Painlev\'e equations: their iteration is given by translations on lattices generated by affine reflection groups.

 \section{Initial-Value Spaces}\label{s:ivs}
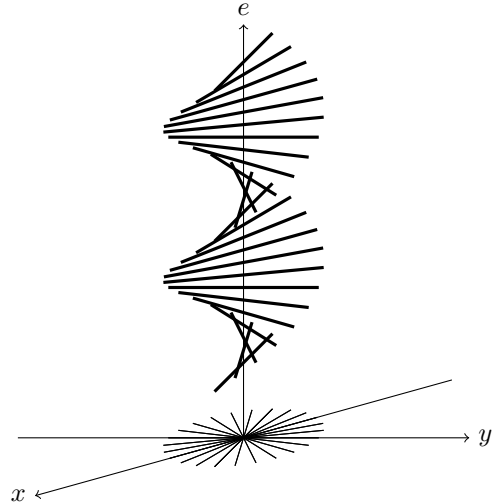
\begin{figure}[H]
    \centering
  \begin{tikzpicture}[scale=0.5]
xc
    \draw[color=black, thin, ->] (0,-2,0) -- (0,9,0) node[above] {$e$}; 
    \draw[color=black, thin, ->] (-6,-2,0) -- (6,-2,0) node[right] {$y$};
    \draw[color=black, thin, ->]  (4,-2,-4) -- (-4,-2,4) node[left] {$x$}; 

    \foreach \i in {0,15,...,360}
    {
        \draw[color=black,thin]  ({-1*2*sin(\i)},-2,{-1*2*cos(\i)}) -- ({2*sin(\i)},-2,{2*cos(\i)});
        \draw[color=black,very thick]  ({-1*2*sin(\i)},{\i / 45},{-1*2*cos(\i)}) -- ({2*sin(\i)},{\i / 45},{2*cos(\i)});
      }
   \end{tikzpicture}
   \caption{The lines in the horizontal $(x, y)$-plane denote local solution trajectories all intersecting at the origin. The vertical axis represents the exceptional line $e$, which replaces the origin after a resolution (i.e., after the transformation in Equation \eqref{eq:2dblowupatorigin} is applied). The bold lines intersecting with the vertical axis represent the resolution of the original solution trajectories.}\label{fig:blowup}
 \end{figure}

 Solutions of dynamical systems follow curves in the space of all possible initial values. But the continuation of a given trajectory is not always guaranteed. In particular, trajectories can all intersect at problematic points called {\em base points} (shown on the horizontal plane in Figure \ref{fig:blowup}). We explain how to disentangle them by using a procedure called {\em blow-up} or resolution in algebraic geometry.
\noindent 

Metaphorically, the solution trajectories are like roads that have a gridlock at the origin. Cars travelling on each road are stuck at the base point because their equations of motion are not defined at the intersection. The mathematical resolution of the gridlock acts by lifting each road to a different vertical height on $e$, before returning the car to the continuation of the road it was on originally. 

  Consider the case of trajectories on the horizontal plane in Figure \ref{fig:blowup}. This plane represents $\mathbb R^2$ containing a one-parameter family of lines intersecting at the origin. A blow up of the origin $(x, y)=(0,0)$ is achieved by taking new coordinates:
\begin{equation}\label{eq:2dblowupatorigin}
  \begin{cases}
    x_1&=\displaystyle \frac{x}{y}\,,\\
    y_1&=y,
    \end{cases}\qquad
  \begin{cases}
    x_2&=x,\\
    y_2&=\displaystyle \frac{y}{x}\,.
    \end{cases}
  \end{equation}
  The line  $e=\{y_1=0\}\cup\{x_2=0\}$ plays a distinguished role and is called an \emph{exceptional line}. It is given by the vertical axis in Figure \ref{fig:blowup}.

Take, for example, an autonomous version of the first Painlev\'e equation, for a given constant $g_2$, in the form $y''=6y^2-g_2/2$. Integration yields
\begin{equation}\label{eq:wpencil}
f(x, y)=(y')^2-4 y^3-g_2 y - g_3,
  \end{equation}
  where $g_3$ is an arbitrary constant. The zeroes of $f$ form a Weierstrass cubic pencil, well known in the theory of elliptic curves \cite{cassels1991lmsst}. In homogeneous coordinates $[u:v:w]\in \mathbb P^2$, we have the homogeneous pencil
  \[
F[u, v, w; g_3]=wv^2-4u^3-g_2uw^2-g_3w^3=0,
\]
where $y=u/w$, $y'=v/w$ are affine coordinates in the finite domain. All curves in the pencil intersect at $[u,v,w]=[0,1,0]$, which lies at infinity. Note that this is precisely where the solution $y(x)$ has a pole.

For the Weierstrass pencil, nine blow-ups are needed to regularize the space. Let the sequence of blow-ups be $\pi_i:X_i\to X_{i-1}$ of $p_i\in X_{i-1}$, with $X=X_9\to X_8\to \ldots\to X_0=\mathbb P^2$. 
Each blow-up $\pi_i$ replaces a base point $p_i$ with an exceptional line $L_i$. At the end of the resolution procedure, we are left with divisor classes $ L_0, \ldots,  L_9$, which form a free $\mathbb Z$-module basis of the Picard group Pic$(X)$, equipped with an intersection form.

A parallel construction can be carried out for the Painlev\'e and discrete Painlev\'e equations. Each point in the regularized space is given by initial values of the equation at some time: $(w, w')$ for the case of ODEs or $(\underline w, w)$ for difference equations. For this reason,  Okamoto \cite{okamoto1979} called the resulting regularized space $X$ for the Painlev\'e equations a \emph{space of initial values}.

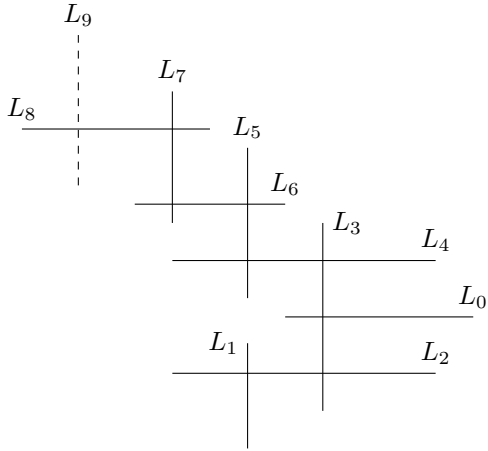
\begin{figure}[H]
    \centering
  \begin{tikzpicture}[scale=0.5]
    \draw[color=black] (-2,0) -- (5,0) node[above] {$ L_2$}; 
    \draw[color=black] (0,-2) -- (0,0.8) node[left] {$ L_1$};
    \draw[color=black]  (2,-1) -- (2,4) node[right] {$ L_3$}; 
    \draw[color=black]  (1,1.5) -- (6,1.5) node[above] {$ L_0$}; 
    \draw[color=black]  (-2,3) -- (5,3) node[above] {$ L_4$}; 
    \draw[color=black] (0,2) -- (0,6) node[above] {$ L_5$};
    \draw[color=black]  (-3,4.5) -- (1,4.5) node[above] {$ L_6$}; 
    \draw[color=black]  (-2,4) -- (-2,7.5) node[above] {$ L_7$}; 
    \draw[color=black]  (-1,6.5) -- (-6,6.5) node[above] {$ L_8$}; 
    \draw[color=black, dashed]  (-4.5,5) -- (-4.5,9) node[above] {$ L_9$}; 
   \end{tikzpicture}
   \caption{The intersection diagram of exceptional divisors arising from resolutions of $\Pone$. The solid lines $L_0, \ldots, L_8$ are lines of self-intersection $-2$, while the dashed line $L_9$ is a line of self-intersection $-1$.}\label{fig:idE8}
 \end{figure}

\begin{figure}[H]
   \centering
  \begin{tikzpicture}[scale=1.0]
      \draw[fill=black] (0,0) circle (0.8ex) node[below=4pt] {$L_1$};
      \draw[fill=black] (1,0) circle (0.8ex) node[below=4pt] {$L_{2}$};
      \draw[fill=black] (2,0) circle (0.8ex) node[below=4pt] {$L_{3}$};
      \draw[fill=black] (3,0) circle (0.8ex) node[below=4pt] {$L_{4}$};
      \draw[fill=black] (4,0) circle (0.8ex) node[below=4pt] {$L_{5}$};
      \draw[fill=black] (5,0) circle (0.8ex) node[below=4pt] {$L_{6}$};
      \draw[thick] (- 0.1,0) -- ( 1.1 ,0);
      \draw[thick] (0.9,0) -- (2.1,0);
      \draw[thick] (1.9,0) -- (3.1,0);
      \draw[thick] (2.9,0) -- (4.1,0);
      \draw[thick] (3.9,0) -- (5.1,0);
      \draw[thick] (4.9,0) -- (6.1,0);
    \draw[fill=black]  (6,0) circle (0.8ex) node[below=4pt] {$L_{7}$};
    \draw[fill=black]  (2,1) circle (0.8ex) node[above=4pt] {$L_{0}$};
    \draw[thick] (2 , 0) -- (2,1);
    \draw[thick] (6 , 0) -- (7.1, 0);
    \draw[fill=black]  (7,0) circle (0.8ex) node[below=4pt] {$L_{8}$};
  \end{tikzpicture}
   \caption{The Dynkin diagram of $E_8^{(1)}$.}\label{fig:ddE8}
 \end{figure}
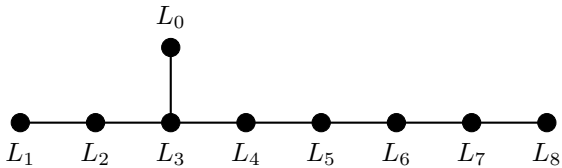

The case of $\Pone$ leads to exceptional lines $ L_0, \ldots,  L_9$ that intersect pairwise as shown in Figure \ref{fig:idE8}. Under the Duval or McKay correspondence, each $L_0, \ldots, L_8$ in Figure \ref{fig:idE8} is mapped to a node, with an edge joining a pair of nodes if the corresponding lines intersect. The resulting graph is the Dynkin diagram of $E_8^{(1)}$ shown in Figure \ref{fig:ddE8}.
\section{Ladder of transcendentality}\label{s:ladder}
In this section, we describe some of the arguments that led to a proof that solutions of Painlev\'e or discrete Painlev\'e equations are indeed new transcendental functions, leading to their position at the top of a ladder of transcendentality. 

Earlier known transcendental functions that are used widely in applied mathematics include the exponential, Airy, Bessel, parabolic cylinder, hypergeometric, or  elliptic functions. These are called {\em classical special functions}. Elliptic functions also satisfy addition formulas (see Figure \ref{f:iter}(\subref{f:iter3})), which are nonlinear difference equations, while the Gamma function solves a linear difference equation. All equations satisfied by such functions are polynomial in the dependent variable and its derivatives or differences.

Such examples laid the foundations of the theory of differential or difference polynomials. Consider any function $y\in C^n$ for any positive integer $n>2$. An $(n+2)$-variable polynomial
\[ P(x, y, y', \ldots , y^{(n)})\,,\] where $y^{(k)}$ is the $k$-th derivative of $y$, is called a  \emph{differential polynomial}. The solutions of $P=0$ are said to be \emph{differentiably algebraic}. There is a natural extension to difference polynomials \[P(n, y_{n}, y_{n+1}, \ldots, y_{n+k})\,,\] where $y_{n+j}$ is the $j$-th iteration of $y_n$.

For example, $y:x\mapsto e^x$ is a root of the differential polynomial $P_0(x, y, y')=y'-y$ and $y:z\mapsto \Gamma(z)$ is a root of the difference polynomial $P_\Gamma\bigl(z, y, \sigma(y)\bigr)=\sigma(y) - z y$. Weierstrass' elliptic function $y=\wp(x)$, with $y'=\wp'(x)$, is a root of the differential polynomial $f$ given in Equation \eqref{eq:wpencil}.

Various approaches have been developed to extend Galois theory, used for number systems, to these settings.  The technical nature of these developments put them outside the scope of this short article, but there are simpler heuristic arguments to show why general solutions of discrete Painlev\'e equations or Painlev\'e equations cannot be rational, algebraic or depend algebraically on classical transcendental functions. 

We outline such an argument for solutions of the first Painlev\'e equation. Consider solutions $w(t)$ of $\Pone$ as functions of initial values. For any nonzero constant $a$, replace $w$ by $u(z)=a^2w(t)$, where $z=a t$, to get
\[
u_{zz}=6 u^2 + a^5 z\,.
\]
If $w(t)$ were rational in its initial values, then $u$ would also be rational. But, in the limit $a\to0$, the solution becomes $\wp(t-t_0, 0, g_3)$, for arbitrary $t_0$, which is not rational in $t_0$ and $g_3$. Therefore, $w(t)$ cannot be a rational function of its initial values.

For the remainder of the argument, we need the fact that solutions of $\Pone$ are meromorphic functions. Consider the possibility that $\Pone$ admits a first integral, which is a polynomial in $w$ and $w'$. In other words, we have a polynomial
\begin{equation}\label{eq:p}
P={w'}^n+Q_1(t, w){w'}^{n-1}+\ldots + Q_n(t, w)=0\,,
\end{equation}
where $Q_i$, with $i=1,\ldots, n$, are polynomial in $w$. This equation should hold in a neighbourhood of each pole $t_0$ of $w$, and transforming to $w=U^{-1}$, $w'=-U'/U^2=-V/U^2$, in a sufficiently small such domain, we find that the equation becomes
\[
  \begin{split}
    V^n-U^2Q_1(t, U^{-1})V^{n-1}+\ldots \\
    + (-1)^nQ_n(t, U^{-1})U^{2n}=0\,.
    \end{split}
\]
The coefficients should be polynomial in $U$. This means that ${\rm deg}_wQ_j(t, w)$ should be at most $2j$.

On the other hand, replacing $w$ and $w'$ in Equation \eqref{eq:p} by $w=a^{-2} u(z)$, $w'=a^{-3} u_z(z)$, where $z=z_0+a t$ and $a$ is a constant, we obtain an equation of the form
\[
P=a^{-k}P_0(u, u_z)+\mathcal O(a^{-k+1}),
\]
for some integer $k$, where $P_0$ is a first integral of the equation $u_{zz}=6 u^2$. Therefore, it must have the form
\[
P_0=b\bigl({u_z}^2-4 u^3\bigr)^d,
\]
where $b$ and $d$ are constants. Putting this together with the above scaling, we find that $k=6d=3n$.

The argument is completed by using the fact that the solutions of $\Pone$ have movable double poles. Substituting the Laurent expansions of $w$ and $w'$ into $P$, we find a contradiction to the requirement that $P$ be polynomial in $w$. Therefore, $P$ cannot exist in the form assumed.

Now assume that $w$ is an algebraic function of $t$. Then it can be expanded in a Puiseux series
\[
w \sim \sum_{n=0}^\infty a_n t^{\rho-n},
\]
for some sector in the domain $|t|\gg 1$. If $\rho\le 0$, then $w$, $w'$, $w''$ would be finite at $z=\infty$, which contradicts the fact that $w$ solves $\Pone$. On the other hand, if $\rho>0$, then it follows from the equation $\Pone$ that it must be integer, but the term $w^2$ in the equation then introduces a larger term of order $\mathcal O(z^{2\rho})$, which is not balanced by any other term. So the equation cannot be satisfied and, therefore, $w$ cannot be algebraic.

Finally, consider the possibility that $w$ is a classical transcendent, which satisfies an algebraic differential equation. Then by further differentiation if necessary, we can eliminate $w''$ and any higher derivatives to obtain a polynomial equation
\[
P(t, w, w')=0.
\]
But this was shown to be impossible above. Therefore, $w$ cannot be a known transcendental function.

\end{document}